\documentclass{amsart}

   \usepackage[
    top    = 2.0cm,
    bottom = 1.95cm,
    left   = 2cm,
    right  = 2cm]    {geometry}

\usepackage{amssymb}
\usepackage{stmaryrd}
\usepackage{fullpage}
\usepackage[colorlinks]{hyperref}

\newtheorem{theorem}{Theorem}[section]

\theoremstyle{definition}

\theoremstyle{remark}

\newtheorem{tab}[theorem]{\bf Table}

\numberwithin{equation}{section}

\newcommand{\FF}{{\mathbb{F}}}

\newcommand{\bC}{{\mathbf{C}}}

\newcommand{\bZ}{{\mathbf{Z}}}

\newcommand{\Aut}{{\operatorname{Aut}}}

\newcommand{\Sp}{{\operatorname{Sp}}}

\newcommand{\GL}{{\operatorname{GL}}}
\newcommand{\GO}{{\operatorname{GO}}}
\newcommand{\Gal}{{\operatorname{Gal}}}

\newcommand{\GF}{\mbox{GF}}

\newcommand{\Magma}{{\sc Magma}}

\begin{document}

\title{Regular orbits of finite primitive solvable groups, the final classification}

\author{Derek Holt}
\address{Mathematics Institute, University of Warwick, Coventry CV47AL, UK.}

\makeatletter
\email{DerekHol127@gmail.com}

\author{YONG YANG}
\address{Department of Mathematics, Texas State University, 601 University Drive, San Marcos, TX 78666, USA.}

\makeatletter
\email{yang@txstate.edu}
\makeatother

\Large
\date{}




\begin{abstract}

    Suppose that a finite solvable group $G$ acts faithfully, irreducibly and quasi-primitively on a finite vector space $V$, and $G$ is not metacyclic.  Then $G$ always has a regular orbit on $V$ except for a few  ``small" cases. We completely classify these cases in this paper.

\end{abstract}

\maketitle
\Large
\section{Introduction} \label{sec:introduction}

Let $G$ be a finite group and $V$ a finite, faithful and completely reducible $G$-module. It is a classical theme to study orbit structure of $G$ acting on $V$. One of the most important and natural questions about orbit structure is to establish the existence of an orbit of a certain size. For a long time, there has been a deep interest and need to examine the size of the largest possible orbits in linear group actions. The orbit $\{v^g \ |\ g \in G\}$ is called regular, if $\bC_G(v)=1$ holds or equivalently the size of the orbit $v^G$ is $|G|$.
The existence of regular orbits has been studied extensively in the literature with many applications to some important questions of character theory and conjugacy classes of finite groups.

In ~\cite{PalfyPyber}, P\'alfy and Pyber asked if it is possible to classify all pairs $A$, $G$ with $(|A|, |G|)=1$ such that $A \leq \Aut(G)$ has a regular orbit on $G$. While the task is pretty challenging, at least for primitive solvable linear groups, we can say something along this line.

Suppose that a finite solvable group $G$ acts faithfully, irreducibly and quasi-primitively on a finite vector space $V$ of dimension $n$ over a finite field of order $q$ and characteristic $p$. Then $G$ has a uniquely determined normal subgroup $E$ which is a direct product of extraspecial $p$-groups for various $p$. We denote $e=\sqrt{|E/\bZ(E)|}$ (an invariant measuring the complexity of the group). It is proved in \cite[Theorem 3.1]{YY2} and \cite[Theorem 3.1]{YY3} that if $e=5,6,7$ or $e \geq 10$ and $e \neq 16$, then $G$ always has regular orbits on $V$. The information on the existence of a regular orbit has been used by several authors to study a variety of problems in the field (for example  ~\cite{DOLFIEmanuele1,DOLFIEmanuele2,LNW,Ponomarenko,YY4,YY5}).

If $e=1$, then $E$ is trivial and $G\leq \Gamma(q^n)$, and it is possible that $G=\Gamma(q^n)$, while $\Gamma(q^n)$ does not have a regular orbit for $n\geq 2$. So for $e=1$ one cannot expect that $G$ necessarily possesses a regular orbit. In this case $G$ is metacyclic and thus there are infinitely many metacyclic primitive linear groups that do not have regular orbits.

There are also other examples for $e>1$, when $G$ does not possess a regular orbit. In ~\cite{YY7}, some more detailed calculation further nails down the cases when $e$ is small. The main result of ~\cite{YY7} implies that there are only finite number of cases left.


Note that we know only a few examples of maximal irreducible primitive solvable subgroups of $\GL(V)$ that are not metacyclic and do not possess a regular orbit. In ~\cite{YY7}, Yang and the others provide a much narrower list of possible groups without regular orbit in ~\cite[Table 3.4]{YY7}.

In this paper, with the help of the computer algebra system \Magma ~\cite{Magma}, we are able to obtain a complete classification for these remaining cases.




\section{Notation and Preliminary Results} \label{sec:Notation and Lemmas}





    If $V$ is a finite vector space of dimension $n$ over $\GF(q)$, where $q$ is a prime power, we denote by $\Gamma(q^n)=\Gamma(V)$ the semilinear group of $V$, i.e.,
\[\Gamma(q^n)=\{x \mapsto ax^{\sigma}\ |\ x \in \GF(q^n), a \in \GF(q^n)^{\times}, \sigma \in \Gal(\GF(q^n)/\GF(q))\}.\]


    For structure of a finite solvable group $G$ that acts faithfully, irreducibly and quasi-primitively on an $d$-dimensional finite vector space $V$ over a finite field $\FF$ of characteristic $p$. Please see ~\cite[Theorem 2.1 and Theorem 2.2]{YY7}.

The following two tables are obtained from ~\cite[Table 3.4]{YY7} for convenience, where $b:=\frac{d}{ea}$.

\begin{tab}{Quasi-primitive solvable groups with $d=ea$ which might not have a regular orbit}\label{exception}
\vspace*{3pt}

    \centering
    \begin{tabular}{|c|c|c|c|c|}
        \hline
        No. & $e$ & $p$ & $d$ & $a$ \\
        \hline
        1   & 16  & 3   & 16  & 1   \\
        2   & 16  & 5   & 16  & 1   \\ \hline
        3   & 9   & 2   & 18  & 2   \\
        4   & 9   & 7   & 9   & 1   \\
        5   & 9   & 13  & 9   & 1   \\
        6   & 9   & 2   & 36  & 4   \\
        7   & 9   & 19  & 9   & 1   \\
        8   & 9   & 5   & 18  & 2   \\ \hline
        9   & 8   & 3   & 8   & 1   \\
        10  & 8   & 5   & 8   & 1   \\
        11  & 8   & 7   & 8   & 1   \\
        12  & 8   & 3   & 16  & 2   \\
        13  & 8   & 11  & 8   & 1   \\
        14  & 8   & 13  & 8   & 1   \\
        15  & 8   & 17  & 8   & 1   \\
        16  & 8   & 19  & 8   & 1   \\
        17  & 8   & 5   & 16  & 2   \\
        18  & 8   & 3   & 24  & 3   \\ \hline
        19  & 4   & 3   & 4   & 1   \\
        20  & 4   & 5   & 4   & 1   \\
        21  & 4   & 7   & 4   & 1   \\
        22  & 4   & 3   & 8   & 2   \\
        23  & 4   & 11  & 4   & 1   \\
        24  & 4   & 13  & 4   & 1   \\
        25  & 4   & 17  & 4   & 1   \\
        26  & 4   & 19  & 4   & 1   \\
        \hline
    \end{tabular}
    \quad
    \begin{tabular}{|c|c|c|c|c|}
        \hline
        No. & $e$ & $p$ & $d$ & $a$ \\
        \hline
        27  & 4   & 23  & 4   & 1   \\
        28  & 4   & 5   & 8   & 2   \\
        29  & 4   & 3   & 12  & 3   \\
        30  & 4   & 29  & 4   & 1   \\
        31  & 4   & 31  & 4   & 1   \\
        32  & 4   & 37  & 4   & 1   \\
        33  & 4   & 41  & 4   & 1   \\
        34  & 4   & 43  & 4   & 1   \\
        35  & 4   & 47  & 4   & 1   \\
        36  & 4   & 7   & 8   & 2   \\
        37  & 4   & 53  & 4   & 1   \\
        38  & 4   & 59  & 4   & 1   \\
        39  & 4   & 61  & 4   & 1   \\
        40  & 4   & 67  & 4   & 1   \\
        41  & 4   & 71  & 4   & 1   \\
        42  & 4   & 73  & 4   & 1   \\
        43  & 4   & 3   & 16  & 4   \\
        44  & 4   & 11  & 8   & 2   \\
        45  & 4   & 5   & 12  & 3   \\
        46  & 4   & 13  & 8   & 2   \\
        47  & 4   & 3   & 20  & 5   \\ \hline
        48  & 3   & 2   & 6   & 2   \\
        49  & 3   & 7   & 3   & 1   \\
        50  & 3   & 13  & 3   & 1   \\
        51  & 3   & 2   & 12  & 4   \\
        52  & 3   & 19  & 3   & 1   \\
        \hline
    \end{tabular}
    \quad
    \begin{tabular}{|c|c|c|c|c|}
        \hline
        No. & $e$ & $p$ & $d$ & $a$ \\
        \hline
        53  & 3   & 5   & 6   & 2   \\
        54  & 3   & 7   & 6   & 2   \\
        55  & 3   & 2   & 18  & 6   \\
        56  & 3   & 11  & 6   & 2   \\
        57  & 3   & 13  & 6   & 2   \\
        58  & 3   & 2   & 24  & 8   \\
        59  & 3   & 17  & 6   & 2   \\
        60  & 3   & 7   & 9   & 3   \\
        61  & 3   & 19  & 6   & 2   \\ \hline
        62  & 2   & 3   & 2   & 1   \\
        63  & 2   & 5   & 2   & 1   \\
        64  & 2   & 7   & 2   & 1   \\
        65  & 2   & 3   & 4   & 2   \\
        66  & 2   & 11  & 2   & 1   \\
        67  & 2   & 13  & 2   & 1   \\
        68  & 2   & 17  & 2   & 1   \\
        69  & 2   & 19  & 2   & 1   \\
        70  & 2   & 23  & 2   & 1   \\
        71  & 2   & 5   & 4   & 2   \\
        72  & 2   & 3   & 6   & 3   \\
        73  & 2   & 29  & 2   & 1   \\
        74  & 2   & 7   & 4   & 2   \\
        75  & 2   & 3   & 8   & 4   \\
        76  & 2   & 11  & 4   & 2   \\
        77  & 2   & 5   & 6   & 3   \\
        78  & 2   & 13  & 4   & 2   \\
        \hline
    \end{tabular}
    \quad
    \begin{tabular}{|c|c|c|c|c|}
        \hline
        No. & $e$ & $p$ & $d$ & $a$ \\
        \hline
        79  & 2   & 3   & 10  & 5   \\
        80  & 2   & 17  & 4   & 2   \\
        81  & 2   & 7   & 6   & 3   \\
        82  & 2   & 19  & 4   & 2   \\
        83  & 2   & 23  & 4   & 2   \\
        84  & 2   & 5   & 8   & 4   \\
        85  & 2   & 3   & 12  & 6   \\
        86  & 2   & 29  & 4   & 2   \\
        87  & 2   & 31  & 4   & 2   \\
        88  & 2   & 11  & 6   & 3   \\
        89  & 2   & 37  & 4   & 2   \\
        90  & 2   & 41  & 4   & 2   \\
        91  & 2   & 43  & 4   & 2   \\
        92  & 2   & 3   & 14  & 7   \\
        93  & 2   & 13  & 6   & 3   \\
        94  & 2   & 47  & 4   & 2   \\
        95  & 2   & 7   & 8   & 4   \\
        96  & 2   & 53  & 4   & 2   \\
        97  & 2   & 5   & 10  & 5   \\
        98  & 2   & 59  & 4   & 2   \\
        99  & 2   & 61  & 4   & 2   \\
        100 & 2   & 67  & 4   & 2   \\
        101 & 2   & 17  & 6   & 3   \\
        102 & 2   & 71  & 4   & 2   \\
        103 & 2   & 73  & 4   & 2   \\
        \hline
    \end{tabular}
\end{tab}

\begin{tab}{Quasi-primitive solvable groups with $b = \frac{d}{ea} > 1$ which might not have a regular orbit}\label{exception2}
\vspace*{3pt}

    \centering
    \begin{tabular}{|c|c|c|c|c|c|}
        \hline
        No. & $e$ & $p$ & $d$ & $a$ & $b$\\
        \hline
         104 & 2 & 3 & 4 & 1 & 2\\
         105 & 2 & 5 & 4 & 1 & 2\\
         106 & 2 & 7 & 4 & 1 & 2\\
         107 & 2 & 11 & 4 & 1 & 2\\
         108 & 2 & 13 & 4 & 1 & 2\\
         109 & 2 & 17 & 4 & 1 & 2\\
         110 & 2 & 3 & 8 & 2 & 2\\
         111 & 2 & 3 & 6 & 1 & 3\\
         112 & 2 & 5 & 6 & 1 & 3\\
         113 & 2 & 3 & 8 & 1 & 4\\
        \hline
     \end{tabular}
    \quad
    \begin{tabular}{|c|c|c|c|c|c|}
        \hline
        No. & $e$ & $p$ & $d$ & $a$ & $b$\\
        \hline
         114 & 3 & 7 & 6 & 1 & 2\\
         115 & 3 & 2 & 12 & 2 & 2\\
         116 & 3 & 2 & 18 & 2 & 3\\
        \hline
         117 & 4 & 3 & 8 & 1 & 2\\
         118 & 4 & 5 & 8 & 1 & 2\\
         119 & 4 & 7 & 8 & 1 & 2\\
         120 & 4 & 11 & 8 & 1 & 2\\
         121 & 4 & 3 & 12 & 1 & 3\\
         122 & 4 & 3 & 16 & 1 & 4\\
         123 & 4 & 3 & 16 & 2 & 2\\
        \hline
     \end{tabular}
    \quad
    \begin{tabular}{|c|c|c|c|c|c|}
        \hline
        No. & $e$ & $p$ & $d$ & $a$ & $b$\\
        \hline
         124 & 8 & 3 & 16 & 1 & 2\\
         125 & 8 & 5 & 16 & 1 & 2\\
         126 & 8 & 3 & 24 & 1 & 3\\
        \hline
         127 & 9 & 2 & 36 & 2 & 2\\
        \hline
    \end{tabular}
\end{tab}

\section{Computations}
In this section, we describe how we constructed candidates for groups $G$ with
parameters equal to one of the entries in Tables~{\em \ref{exception}}
and~{\em \ref{exception2}} on a computer, and checked in each case whether
there were any such examples without regular orbits. We carried out these
computations in \Magma. The results of these computations are tabulated in
the next section.

We know from ~\cite[Theorem 2.2]{YY7} that $G$ has a normal subgroup $F$,
which is a central product of a subgroup $U=Z(F)$ and an extraspecial group
$E$ of order $r^{2e+1}$, where $|U|$ divides $p^a-1$ and $U$ acts irreducibly
on a subspace $W$ of $V$ of dimension $a$. Since $G$ is quasi-primitive,
$U$ acts homogeneously on $V$, and by \cite[Lemma 1.10]{LucMenMor01}
(applied with $M$, $K$ and $F$ equal to $U$, $\FF_p$ and $\FF_{p^a}$),
we can regard $V$ as a vector space over the field $\FF_{p^a}$ of order $p^a$,
and we have $A = C_G(U) \le C_{\GL(d,p)}(U) \cong \GL(d/a,p^a)$.

Furthermore, $G$ is isomorphic to a subgroup of the normalizer of
$C_{\GL(d,p)}(U)$ in $\GL(d,p)$, which is isomorphic to $\Gamma L(d/a,p^a)$
(i.e. the split extension $\GL(d/a,p^a):\langle \sigma \rangle$, where
$\sigma$ acts on $\GL(d/a,p^a)$ as a field automorphism of order $a$)
and $G/A$ can be identified with a subgroup of $\Gal(\FF_{p^a},\FF_p)$,
which is cyclic of order $a$.

We shall now summarize some properties of extraspecial and symplectic-type
groups and their representations.  Convenient background references for much
of this material are \cite[Section 4.6]{KL90} or \cite[Section 5.5]{Gorenstein}.

For a prime $r$ and integer $e \ge 1$ there are two isomorphism types of
extraspecial $r$-groups $E$ of order $r^{2e+1}$, and they both arise as central
products of $e$ extraspecial groups of order $r^3$.
Their faithful absolutely irreducible representations in characteristics other
than $r$ have dimension $e$ and are quasi-equivalent to each other
(i.e. equivalent under the action of $\Aut(E)$).

Assume first that $ea=d$ (or, equivalently, that $b=1$). Then, since $e$
is the dimension of all non-linear absolutely irreducible representations of
$E$ in characteristic $p \ne r$, the group $E$ must be absolutely
irreducible as a subgroup of $\GL(e,p^a)$ and, for a given isomorphism type
of $E$, since its faithful absolutely irreducible representations are
quasi-equivalent, there is a unique conjugacy class of subgroups of
$\GL(e,p^a)$ isomorphic to $E$.

Our methods for the cases when $r$ is odd and even are slightly different,
so we consider them separately.
Suppose first that $r$ is odd.  We claim that $E$ must have exponent $r$.
The other isomorphism type of extraspecial group has exponent
$r^2$, and its elements of order $r$ form a characteristic subgroup $E_r$ of
index $r$ in $E$ with non-cyclic center of order $r^2$.  So $E_r$ has no
faithful irreducible irreducible representations, but it acts faithfully
on $V$, so it cannot be acting homogeneously, contradicting the
quasi-primitivity of $G$.

There is existing functionality in \Magma\ for constructing $E$ as a subgroup
of $\GL(e,p^a)$ and its normalizer $N_A$ in $\GL(e,p^a)$ (which is not
usually a solvable group). The group $N_A$ has the structure
$Z_0 r^{1+2e}.\Sp(2e,r)$ (with $|E| = r^{1+2e}$), where $Z_0:= Z(\GL(d/a,p^a))$
is the group of scalar matrices. The group $E$ consists of the elements of
order dividing $r$ in $O_r(N_A)$, and so it is characteristic in $N_A$.

After constructing $N_A$, we embed it in $\GL(d,p)$ using the natural
embedding $\GL(e,p^a) \to \GL(d,p)$. Then, as a subgroup of  $\GL(d,p)$,
$Z_0E$ acts irreducibly with centralizing field $\FF_{p^a}$, so
$C_{\GL(d,p)}(Z_0E) = Z_0$.  Note that the normal subgroup $A$ of the
group $G$ that we are attempting to construct is the intersection of $G$ with
$N_A$, so $G$ is a subgroup of $N := N_{\GL(d,p)}(N_A)$, and the method that
we chose to find $G$ involves computing this group $N$. To do that,
we compute $\Aut(N_A)$, and then check which outer automorphisms of $N_A$ can
be induced by conjugation in $\GL(d,p)$. (This uses the fact that
$C_N(N_A) \le N_A$, which follows from $C_{\GL(d,p)}(Z_0E) = Z_0 \le N_A$.)
This automorphism group computation was one of the slowest parts of the
complete process, and it is possible that there are faster ways of
computing $N$ from $N_A$, but it eventually completed successfully in all
of the examples.

After computing $N$, we compute its subgroups of increasingly large index,
by repeated application of the \textsf{MaximalSubgroups} command
in \Magma, using conjugacy testing to ensure that we only consider one
representative of each $N$-conjugacy class of subgroups. For each such
subgroup, we test whether it is solvable and quasi-primitive. If so, then we
test whether it has regular orbits. If so then we do not need to consider any of
its proper subgroups, because they would also have regular orbits. If not,
then we have identified an example without regular orbits.

The situation is more complicated when $r=2$.  In that case, the extraspecial
groups of order $2^3$ are $Q_8$ and $D_8$, the dihedral and quaternion groups
of order $8$, and those of order $2^{2e+1}$ are $E^+$, a central product of
$e$ copies of $D_8$, and $E^-$, a central product of $e-1$ copies of $D_8$
and one of $Q_8$.  (Note that $D_8*D_8 \cong Q_8*Q_8$.)
Their central products $E^+*C_4$ and $E^-*C_4$ with a cyclic group of order
four are isomorphic \emph{symplectic-type} groups $S$ of order $2^{2e+2}$.
The faithful absolutely irreducible representations of $E^+$, $E^-$, and $S$
in characteristic $p \ne 2$ are quasi-equivalent and have dimension $e$.
Those of $E^+$ and $E^-$ can be written over $\GF(p^a)$ for any odd prime
$p$ and any $a \ge 1$, whereas those of $S$ can be written over $\GF(p^a)$
if and only if $4 \mid p^a-1$; i.e. if and only if either $p \equiv 1 \bmod 4$,
or  $p \equiv 3 \bmod 4$ and $a$ is even.

Suppose first that $4 \mid p^a-1$. Then the normalizer $N(S)_A$ of $S$ in
$\GL(e,p^a)$ contains the normalizers of $E^+$ and of $E^-$, and so we deal
with both of these cases together by computing it. The group $N(S)_A$
has the structure $Z_0 2^{1+2e}.\Sp(2e,2)$, and $S$ consists of the elements
of $O_2(N(S)_A)$ of order dividing $4$, so $S$ is characteristic in $N(S)_A$.
We use the same process as for the case with $r$ odd but with $S$ in place of
$E$.

When $4$ does not divide $p^a-1$, we have $E = O_2(G) \lhd G$, where $E$
can be isomorphic to either $E^+$ or $E^-$, and we must carry out the
computations for these two cases separately.
We have $N_A = N_{\GL(e,p^a)}(E) \cong  Z 2^{1+2e}.\GO^+(2e,2)$ and
$N_0 = Z2^{1+2e}.\GO^-(2e,2)$ in the two cases, and we proceed as in
the case $r$ odd in both cases.

It remains to consider the case $d/(ea) = b > 1$. Then, by quasi-primitivity,
the group $E$ acts homogeneously as a subgroup of $\GL(d,p)$ and, since
it is centralized by $U$, which acts as scalar multiplication as a
subgroup of $\GL(d/a,p^a)$, the group $E$ also acts homogeneously as a
subgroup of  $\GL(d/a,p^a)$. So it has $b$ isomorphic absolutely
irreducible constituents, each of dimension $e$ over $\FF_{p^a}$.
Now, by \cite[Theorem 3.5.4]{Gorenstein}, we have
$C_A := C_{\GL(d/a,p^a)}(E) \cong \GL(b,p^a)$, and hence also
$CC_A := C_{\GL(d/a,p^a)}(C_A) \cong \GL(e,p^a)$. Now the normalizer in
$\GL(d,p)$ of $E$ also normalizes $N_{CC_A}(E) \cong N_{\GL(e,p^a)}(E)$,
and hence it also normalizes the subgroup
$N_A := \langle C_A, N_{CC_A}(E) \rangle$
(or $N(S)_A := \langle C_A, N_{CC_A}(S) \rangle$  when $4 \mid p^a-1$).

We can compute $N_{CC_A}(E)$ (or $N_{CC_A}(S)$)  as in the case $b=1$, and
$C_A$ is also straightforward to compute, so we can compute the group $N_A$
(or $N(S)_A$), and we use this to construct the normalizer in
$\GL(d,p)$ of $E$ and its subgroups in the same way as in the case $b=1$.

\section{Table of results}
Here is a list of those entries in Table~{\em \ref{exception}} for which there
is at least one example of a group with no regular orbit. In each case
we give the number ``num gps'' of such examples (up to conjugacy in
$\GL(d,p)$), and the order ``max $|G|$'' of the largest example.
In cases where $r=2$ and $4$ does not divide $2^a-1$, we have handled
the $E^+$ and $E^-$ cases separately. 

The whole data package is provided in a separate file.

\smallskip
\begin{tab}{Parameters of quasi-primitive solvable groups that do not have a regular orbit}\label{exceptiontrue}
\vspace*{3pt}

    \centering
    \begin{tabular}{|c|c|c|c|c|c|c|r|c|}
        \hline
        No. & $e$ & $p$ & $d$ & $a$ & $b$ & num gps & max $|G|$ & Note \\
        \hline
 1 & 16 & 3 & 16 & 1 & 1 & 12 & 15925248 & $E^-$ \\
 3 &  9 & 2 & 18 & 2 & 1 & 40 & 559872 &\\
 9 &  8 & 3 &  8 & 1 & 1 & 27 & 18432 & $E^+$\\
 9 &  8 & 3 &  8 & 1 & 1 & 71 & 165888 & $E^-$\\
10 &  8 & 5 &  8 & 1 & 1 & 22 & 331776 & \\
19 &  4 & 3 &  4 & 1 & 1 & 14 & 2304 & $E^+$\\
19 &  4 & 3 &  4 & 1 & 1 & 9  & 640 & $E^-$  \\
20 &  4 & 5 &  4 & 1 & 1 & 24 & 4608 &  \\
21 &  4 & 7 &  4 & 1 & 1 & 17 & 6912 & $E^+$  \\
22 &  4 & 3 &  8 & 2 & 1 & 72 & 18432 &  \\
23 &  4 & 11 & 4 & 1 & 1 & 4  & 11520 & $E^+$  \\
24 &  4 & 13 & 4 & 1 & 1 & 5  & 13824 &  \\
25 &  4 & 17 & 4 & 1 & 1 & 4  & 18432 &  \\
28 &  4 & 5 &  8 & 2 & 1 & 3  & 55296 &  \\
48 &  3 & 2 &  6 & 2 & 1 & 7  & 1296 &  \\
49 &  3 & 7 &  3 & 1 & 1 & 4  & 1296 &  \\
50 &  3 &13 &  3 & 1 & 1 & 2  & 2592 &  \\
51 &  3 & 2 & 12 & 4 & 1 & 8  & 12960 &  \\
52 &  3 &19 &  3 & 1 & 1 & 1  & 3888 &  \\
53 &  3 & 5 &  6 & 2 & 1 & 10 & 10368 &  \\
62 &  2 & 3 &  2 & 1 & 1 & 2  & 48 &  \\
63 &  2 & 5 &  2 & 1 & 1 & 2  & 96 &  \\
64 &  2 & 7 &  2 & 1 & 1 & 2  & 144 &  \\
65 &  2 & 3 &  4 & 2 & 1 & 13 & 384 &  \\
66 &  2 & 11 &  2 & 1 & 1 & 2 & 240 &  \\
67 &  2 & 13 &  2 & 1 & 1 & 2 & 288 &  \\
68 &  2 & 17 &  2 & 1 & 1 & 3 & 384 &  \\
69 &  2 & 19 &  2 & 1 & 1 & 2 & 432 &  \\
71 &  2 & 5 &  4 & 2 & 1 & 16 & 1152 &  \\
72 &  2 & 3 &  6 & 3 & 1 & 2  & 1872 &  \\
74 &  2 & 7 &  4 & 2 & 1 & 7  & 2304 &  \\
75 &  2 & 3 &  8 & 4 & 1 & 10 & 7680 &  \\
117 & 4 & 3 & 8 & 1 & 2 & 9 & 2304&\\
        \hline
    \end{tabular}
\end{tab}

\section{Acknowledgement} \label{sec:Acknowledgement}
The research of the first author  was partially supported by the NSFC (No: 11671063) and a grant from the Simons Foundation (No 499532).

\small


\begin{thebibliography}{19}
    \bibitem{Magma}
Wieb Bosma, John Cannon, Catherine Playoust.
The {M}agma algebra system. {I}. {T}he user language.
\textit{J. Symbolic Comput.}, 24(3-4) (1997), 235--265.
\newblock Computational algebra and number theory (London, 1993).
    
    \bibitem {DOLFIEmanuele1} {S. Dolfi and E. Pacifici}, `Zeros of Brauer characters and linear actions of finite groups', \textit{J. Algebra} 340 (2011), 104--113.
    \bibitem {DOLFIEmanuele2} {S. Dolfi and E. Pacifici}, `Zeros of Brauer characters and linear actions of finite groups: small primes', \textit{J. Algebra} 399 (2014), 343--357.


\bibitem{Gorenstein} D. Gorenstein, \emph{Finite Groups}, Harper and Row, 1968.


\bibitem{KL90}{P. Kleidman and M. Liebeck}, \emph{The subgroup  structure of the finite classical groups},
Cambridge University Press, Cambridge, 1990.

    \bibitem{LNW} {M. Lewis, G. Navarro, T.\,R. Wolf}, `$p$-parts of character degrees and the index of the Fitting subgroup', \textit{J. Algebra} 411 (2014), 182--190.

\bibitem{LucMenMor01} {A. Lucchini, F. Menegazzo, M. Morigi}, `On the number of generators and composition length of finite linear groups', \textit{J. Algebra} 243 (2001), 427--447.



    \bibitem{PalfyPyber} {P. P\'alfy and L. Pyber}, `Small groups of automorphisms', \textit{Bull. Lond. Math. Soc.} 30 (1998), 386--390.

    \bibitem {Ponomarenko} {I.\,N. Ponomarenko}, `Bases of Schurian antisymmetric coherent configurations and isomorphism test for Schurian tournaments',  \textit{J. Math. Sci. (N.Y.)} 192 (2013), 316--338.



    \bibitem{YY2} {Y. Yang}, `Regular orbits of finite primitive solvable groups', \textit{J. Algebra} 323 (2010), 2735--2755.

    \bibitem{YY3} {Y. Yang}, `Regular orbits of finite primitive solvable groups, II', \textit{J. Algebra} 341 (2011), 23--34.

    \bibitem{YY4} {Y. Yang}, `Regular orbits of nilpotent subgroups of solvable linear groups', \textit{J. Algebra} 325 (2011), 56--69.

    \bibitem{YY5} {Y. Yang}, `Large character degrees of solvable $3'$-groups', \textit{Proc. Amer. Math. Soc.} 139 (2011), 3171--3173.

    \bibitem{YY6} {Y. Yang}, `Blocks of small defect', \textit{J. Algebra} 429 (2015), 192--212.
    \bibitem{YY7} {Y. Yang, A. Vasil'ev, E. Vdovin}, `Regular orbits of finite primitive solvable groups, III', \textit{J. Algebra} 590 (2022), 139--154.
\end{thebibliography}
\end{document}